\documentclass[sn-basic, Numbered]{sn-jnl}
\usepackage{graphicx}%
\usepackage{multirow}%
\usepackage{amsmath,amssymb,amsfonts}%
\usepackage{amsthm}%
\usepackage{mathrsfs}%
\usepackage[title]{appendix}%
\usepackage{xcolor}%
\usepackage{textcomp}%
\usepackage{manyfoot}%
\usepackage{booktabs}%
\usepackage{algorithm}%
\usepackage{algorithmicx}%
\usepackage{algpseudocode}%
\usepackage{listings}%
\usepackage[vcentermath,enableskew]{youngtab}
\usepackage{colortbl}

\theoremstyle{thmstyleone}%
\newtheorem{theorem}{Theorem}
\newtheorem{proposition}[theorem]{Proposition}%
\newtheorem{lemma}{Lemma}%
\newtheorem{corollary}{Corollary}%

\theoremstyle{thmstyletwo}%

\theoremstyle{thmstylethree}%
\newtheorem{definition}{Definition}%

\renewcommand{\k}{\mathbb{F}}

\begin{document}

\title[Article title]{Schur--Weyl dualities for the rook monoid: an approach via Schur algebras}

\author[1]{\fnm{Carlos A.} \sur{M. André}}\email{caandre@ciencias.ulisboa.pt}

\author*[2]{\fnm{Inês} \sur{Legatheaux Martins}}\email{ilegatheaux@yahoo.co.uk}

\affil[1, 2]{\orgdiv{Departamento de Matemática}, \orgname{Faculdade de Ciências da Universidade de Lisboa}, \orgaddress{\street{Campo Grande, Edifício C6, Piso 2}, \postcode{1749-016}, \city{Lisboa}, \country{Portugal}}}

\abstract{The rook monoid, also known as the symmetric inverse monoid, is the archetypal structure when it comes to extend the principle of symmetry. In this paper, we establish a Schur--Weyl duality between this monoid and an extension of the classical Schur algebra, which we name the extended Schur algebra. We also explain how this relates to Solomon's Schur--Weyl duality between the rook monoid and the general linear group and mention some advantages of our approach.}

\keywords{Schur--Weyl duality, rook monoid, Schur algebras, representation theory of associative algebras, tensor spaces}

\pacs[MSC Classification]{20M30, 20G43, 16G99, 16S50, 20M18, 22E46}

\maketitle

\section{Introduction}\label{intro}

\vspace{0.3cm}

Throughout this article, $\k$ is a field of characteristic zero unless explicitly specified and $V$ is a $d$-dimensional vector space over $\k$. The symmetric group $S_{n}$ acts on the tensor space $\otimes^{n} V$ by place permutations. By fixing a basis of $V$, $GL(V)$ can be identified with the general linear group of all $d\times d$ non-singular matrices with entries in $\k$, herein denoted $G_{d}$. If $V$ is the natural module for the group algebra $\k G_{d}$, then $G_{d}$ acts diagonally on $\otimes^{n} V$. This action commutes with that of $S_{n}$ on $\otimes^{n} V$ by place permutations. In case $\k = \mathbb{C}$, I. Schur \cite{Schur1927} established that each action generates the full centralizer of the other on ${\rm End}_{\, \k}(\otimes^{n} V)$, a result which was made popular by H. Weyl \cite{Weyl}. This seminal example of a double centralizer phenomenon, now known as the {\it classical Schur--Weyl duality}, provides a deep insight on the interactions between the representation theories of $G_{d}$ and $S_{n}$.

Results of R. M. Thrall \cite{Thrall_ModTensorsII_44}, C. de Concini and C. Procesi \cite{DeConciniProcesi}, R. Carter and G. Lusztig \cite{CarterLusztig} and J. A. Green \cite{Green1} show that classical Schur--Weyl duality remains true if $\k$ is an infinite field of any characteristic. More recently, classical Schur--Weyl duality was extended to sufficiently large finite fields by D. Benson and S. Doty \cite{BensonDoty_SWDualFF09}.

These results can be better understood in the context of Schur algebras \cite{Green1}. Implicit in zero characteristic in Schur's Ph.D thesis \cite{SchurThesis}, Schur algebras were defined over arbitrary infinite fields by J. A. Green in a seminal monograph \cite{Green1}. The {\it Schur algebra} $\mathcal{S}_{\, \scriptscriptstyle{\k}}(d,n)$ can be identified with the centralizer algebra ${\rm End}_{\, \mathbb{\k}S_{n}}(\otimes^{n} V)$ of $\k S_{n}$ on $\otimes^{n} V$ with respect to place permutations and the family $\left \{ \mathcal{S}_{\, \scriptscriptstyle{\k}}(d,r) \right \}_{r\geq 0}$ completely determines the polynomial representations of $\k G_{d}$. Moreover, $\mathcal{S}_{\, \scriptscriptstyle{\k}}(d,n)$ is an important example of a cellular algebra \cite{GrahamLehrer}. Thus, the classical Schur--Weyl duality can be stated in terms of these finite-dimensional algebras in a very general setting. 

There are numerous other examples of ``Schur--Weyl dualities". For instance, in characteristic zero, the centralizer algebras associated with the diagonal action of subgroups of $G_{d}$ such as the orthogonal group $O_{d}$ and the symmetric group $S_{d}$ on $\otimes^{n} V$ are, respectively, the {\it Brauer algebra} \cite{Brauer} and the {\it partition algebra} \cite{Martin_PartAlg_0, Martin_PartAlg_1, Martin_PartAlg_2, Jones_Potts} (see also \cite{HalvRamPartAlg}). As before, the translation of these results in the language of Schur algebras and their generalizations has widely expanded our knowledge of the properties of these algebras in the modular case (see, among many others, \cite{DipperDotyHu, DotyHu_Orthogonal, BowmanDotyMartin_SchurWeylPartAlg}).

In 2002, L. Solomon \cite{Solomon} established a Schur--Weyl duality between $G_{d}$ and an important finite inverse monoid. Inverse monoids were introduced in \cite{Wagner52} as a natural generalization of groups to deal with aspects of symmetry which the latter couldn't capture (see \cite{Lawson} for further details on this viewpoint). The archetypal example of such a structure is the {\it symmetric inverse monoid}, also called the {\it rook monoid} \cite{Solomontits}.

For our purposes, the rook monoid $R_{n}$ is the set of all bijective partial maps from $\mathbf{n}=\{1, \ldots , n\}$ to itself under the usual composition of partial functions. It contains $S_{n}$ and it is isomorphic to the monoid under matrix multiplication of all $n\times n$ matrices with at most one entry equal to $1_{\k}$ in each row and in each column and zeros elsewhere. It plays the same r\^ole for inverse monoids that $S_{n}$ does for groups and thus is the archetypal structure when it comes to extend the principle of symmetry.

In his influential article \cite{Solomon}, L. Solomon proved that $R_{n}$ acts on tensors by ``place permutations". More precisely, he showed that, if $\k$ has characteristic zero and $d\geq n$, $\k R_{n}$ acts as the centralizer algebra for the action of $G_{d}$ on $\otimes^{n} U$, where $U=V\oplus U_{0}$ is the direct sum of the natural $d$-dimensional module $V$ and the trivial module $U_{0}$. 

Since its publication, this result proved to be a special case of an important Schur-Weyl duality on tensor spaces for the Hecke algebra analog for $R_{n}$, known as the $q$-rook monoid (see \cite{Halvqrook, Solomonqrook, Paget} and references therein). It also influenced other authors into establishing Schur-Weyl dualities between $R_{n}$ and other finite inverse semigroups \cite{GannaMazorchuk}. Moreover, it led to the investigation of a number of interesting algebras. For instance, the centralizer algebras associated with the restriction of the action of $G_{d}$ on $\otimes^{n} U$ to subgroups such as the orthogonal subgroup $O_{d}$ and the symmetric $S_{d}$ are, respectively, the {\it rook Brauer algebra} \cite{HalvDelMas_RookBrauerRep, MazorchukMartin_PartialBrauerAlg} and the {\it rook partition algebra} \cite{Grood_RookPartitionAlg}.

The main purpose of this article is to show that Solomon's Schur--Weyl duality for $R_{n}$ and $G_{d}$ can be stated in terms of an extension of the classical Schur algebra. 

We achieve this by defining an $\k$-algebra $\mathcal{S}_{\,\scriptscriptstyle{\k}}(d, \mathbf{n})$ which we call the {\it extended Schur algebra} and which satisfies 
$$\mathcal{S}_{\,\scriptscriptstyle{\k}}(d,\mathbf{n})\cong \bigoplus_{r=0}^{n} \mathcal{S}_{\, \scriptscriptstyle{\k}}(d,r).$$
We then prove that $\mathcal{S}_{\, \scriptscriptstyle{\k}}(d, \mathbf{n})$ determines the homogeneous polynomial representations of $\k G_{d}$ of degree at most $n$, a result that holds for arbitrary infinite fields. Finally, we establish a Schur--Weyl duality between $\mathcal{S}_{\, \scriptscriptstyle{\k}}(d, \mathbf{n})$ and $R_{n}$ on $\otimes^{n} U$, when $d\geq n$ and $\k$ has zero characteristic. To show that our viewpoint provides a deeper insight on the representation theory of $R_{n}$ and its interactions to those of general linear and symmetric groups, we also mention some applications of our approach.

This paper is organized as follows. Section 2 begins with a brief overview on (split) semisimple algebras, double centralizer theory and classical Schur--Weyl duality. This is followed by a description of structural aspects of the classical Schur algebra $\mathcal{S}_{\, \scriptscriptstyle{\k}}(d,n)$ and an outline of the representation theory of the rook monoid $R_{n}$.

In Section 3, we view $G_{d}\subseteq G_{d+1}$ under a natural embedding and we explain how the restriction of the diagonal action of $G_{d+1}$ on $\otimes^{n} U$ to $G_{d}$ gives rise to the extended Schur algebra $\mathcal{S}_{\, \scriptscriptstyle{\k}}(d, \mathbf{n})$. After describing this algebra's structure, we prove that the module category for $\mathcal{S}_{\, \scriptscriptstyle{\k}}(d, \mathbf{n})$ is equivalent to the category of homogeneous polynomial $G_{d}$-modules of degree at most $n$. Finally, we establish a Schur--Weyl duality on $\otimes^{n} U$ between $\mathcal{S}_{\, \scriptscriptstyle{\k}}(d, \mathbf{n})$ and $\k R_{n}$. We end by explaining how this result relates to Solomon's Schur-Weyl duality \cite{Solomon} and mentioning some consequences of our approach.

We should note that some of the techniques used herein apply to infinite fields of any characteristic. The fact that our main result relies on the semisimplicity of the monoid algebra of $R_{n}$ has made us decide to work in characteristic zero. However, since we hope to treat the modular case in the near future, we have pointed out all the results in this article that remain true for arbitrary infinite fields.

\vspace{0.3cm}

\section{Preliminaries}\label{sec:prelim}

\vspace{0.3cm}

\subsection{Double centralizer theory and classical Schur--Weyl duality}\label{subsec2a}

\vspace{0.3cm}

Henceforth, the term ``module" refers to a finite-dimensional left module unless explicitly stated otherwise and $\k$ is a field of characteristic zero. Let $\mathcal{A}$ be a finite-dimensional split semisimple algebra over $\k$. By classical Artin-Wedderburn theory \cite[Theorem $3.34$]{CurtisReiner81}, this means that there is an isomorphism of $\k$-algebras
$$\mathcal{A}\cong \bigoplus_{\lambda \in \Lambda} \mathcal{M}_{d_{\lambda}}(\k),$$
for some finite index set $\Lambda$ and positive integers $d_{\lambda}$. For each $\lambda\in \Lambda$, there is, up to isomorphism, one simple $\mathcal{A}$-module $S_{\lambda}$ and $\left \{S_{\lambda} : \lambda \in \Lambda \right \}$ is a complete set of representatives of the isomorphism classes of simple modules of $\mathcal{A}$.

If $M$ is a finite-dimensional $\mathcal{A}$-module, its decomposition into simple $\mathcal{A}$-modules is given by
\begin{equation}
\label{eq:decompirred}
M\cong \bigoplus_{\lambda \in \Lambda} m_{\lambda}S_{\lambda},
\end{equation}
where $m_{\lambda}$ is a non-negative integer called the {\it multiplicity} of $\lambda$ in $M$. We say that $\lambda \in \Lambda$ {\it appears} in $M$ if $M$ contains a submodule isomorphic to $S_{\lambda}$ (that is, if  $m_{\lambda}>0$).

Let $\rho: \mathcal{A} \rightarrow {\rm End}_{\, \k}(M)$ be the representation corresponding to the $\mathcal{A}$-module $M$ with decomposition given by Expression (\ref{eq:decompirred}). The {\it centralizer algebra of} $\mathcal{A}$ on $M$ is the finite-dimensional $\k$-algebra 
$${\rm End}_{\, \mathcal{A}}(M)=\{\phi \in {\rm End}_{\, \k}(M) : \phi \rho(a)=\rho(a) \phi, {\rm for \ all \ } a\in \mathcal{A}\}.$$
The action of ${\rm End}_{\, \mathcal{A}}(M)$ on $M$ defined by $\phi \cdot x = \phi(x)$, with $\phi\in {\rm End}_{\, \mathcal{A}}(M)$ and $x\in M$, turns $M$ into an ${\rm End}_{\, \mathcal{A}}(M)$-module. Since the actions of $\mathcal{A}$ and ${\rm End}_{\, \mathcal{A}}(M)$ on $M$ commute, $\rho$ induces a homomorphism of $\k$-algebras $\rho: \mathcal{A} \rightarrow {\rm End}_{\, {\rm End}_{\, \mathcal{A}}(M)}(M)$. In case $M$ is a right $\mathcal{A}$-module, the corresponding homomorphism is obtained via the representation $\rho : \mathcal{A}^{\circ} \rightarrow {\rm End}_{\k}(M)$ of $\mathcal{A}^{\circ}$ on $M$, where $\mathcal{A}^{\circ}$ is the opposite algebra of $\mathcal{A}$.

The next theorem summarizes the basic dual relationship between the $\k$-algebras $\mathcal{A}$ and ${\rm End}_{\mathcal{A}}(M)$ and can be found (with a different formulation) in \cite[Section $3.$D]{CurtisReiner81}.

\vspace{0.2cm}

\begin{theorem}[Double Centralizer Theorem]
\label{teo:doublecent}
Let $\mathcal{A}$ be a finite-dimensional split semisimple $\k$-algebra and let $\{S_{\lambda} : \lambda \in \Lambda\}$ be a complete set of representatives of the isomorphism classes of simple $\mathcal{A}$-modules with $\dim_{\scriptscriptstyle{\k}}(S_{\lambda})=d_{\lambda}$, for all $\lambda \in \Lambda$. Let $M$ be an $\mathcal{A}$-module such that $\displaystyle{M\cong \bigoplus_{\lambda \in \Lambda^{'}} m_{\lambda} S_{\lambda}}$, with $\Lambda^{'} = \{\lambda \in \Lambda : m_{\lambda} > 0\}$. Then:
\begin{enumerate}
\item [$(a)$] there is an isomorphism of $\k$-algebras $\displaystyle{{\rm End}_{\, \mathcal{A}}(M)\cong \bigoplus_{\lambda\in \Lambda^{'}} \mathcal{M}_{m_{\lambda}}(\k)}$; in particular, ${\rm End}_{\, \mathcal{A}}(M)$ is a finite-dimensional (split) semisimple $\k$-algebra;
\item [$(b)$] as an ${\rm End}_{\, \mathcal{A}}(M)$-module, $\displaystyle{M\cong \bigoplus_{\, \lambda \in \Lambda^{'}} d_{\lambda}E_{\lambda}}$, where, for each $\lambda\in \Lambda^{'}$, $E_{\lambda}$ is a simple ${\rm End}_{\, \mathcal{A}}(M)$-module of dimension $m_{\lambda}$;
\item [$(c)$] if $M$ is a faithful $\mathcal{A}$-module, the corresponding representation $\rho:\mathcal{A} \rightarrow {\rm End}_{\, \k}(M)$ induces an isomorphism of $\k$-algebras $$\mathcal{A} \cong {\rm End}_{\, {\rm End}_{\, \mathcal{A}}(M)}(M),$$ i.e., the actions of $\mathcal{A}$ and ${\rm End}_{\mathcal{A}}(M)$ on $M$ generate the full centralizer of the other in ${\rm End}_{\, \k}(M)$.
\end{enumerate}
\end{theorem}

\vspace{0.2cm}

Schur--Weyl duality is a cornerstone of representation theory that amounts to two double centralizer results which involve the symmetric and general linear groups.

Since $\k$ has characteristic zero, it is known that the group algebra of the symmetric group $\k S_{n}$ is a split semisimple algebra of dimension $n \, !$ (see, for example, \cite[Theorem $5.9.$]{JacobsonVol2}). On the other hand, the isomorphism classes of simple modules for $\k S_{n}$ are indexed by partitions of $n$.

By a {\it partition} of $n$, we mean a sequence $\lambda = (\lambda_{1}, \lambda_{2}, \ldots , \lambda_{t})$ of weakly decreasing non-negative integers ($\lambda_{1}\geq \lambda_{2}\geq \ldots \geq \lambda_{t}\geq 0$) whose sum is equal to $n$. We refer to $\lambda_{1}, \ldots, \lambda_{t}$ as the {\it parts} of $\lambda$ and write $\lambda \vdash n$ to indicate that $\lambda$ is a partition of $n$.

Let $V$ be an $\k$-space of dimension $d$. For reasons that will become apparent later, we view $\otimes^{n} V$ as a right $\k S_{n}$-module on which $S_{n}$ acts by place permutations as
\begin{equation}
\label{eq:placepermutsn}
(v_{1}\otimes \ldots \otimes v_{n}) \cdot \sigma = v_{\sigma(1)}\otimes \ldots \otimes v_{\sigma(n)},
\end{equation}
for all $v_{1}, \ldots , v_{n} \in V$ and $\sigma \in S_{n}$.

Let $G_{d}$ be the general linear group of all $d\times d$ invertible matrices with entries in $\k$. If $V$ has as $\k$-basis $\{e_{1}, \ldots , e_{d}\}$, we consider the {\it natural} left action of $G_{d}$ on $V$, defined on basis elements by
\begin{equation}
\label{eq:GaccaoV}
g \cdot e_{j} = \sum_{i=1}^{d} c_{i,j}(g)e_{i}, \ {\rm for \ } g\in G_{d} \ {\rm and \ } j = 1, \ldots , d,
\end{equation}
where, for all $1\leq i, j \leq d$, $c_{i,j}: G_{d}\rightarrow \k$ is the coordinate function which sends a matrix in $G_{d}$ to its $(i,j)$-th entry. It follows that $G_{d}$ acts on $\otimes^{n} V$ via the {\it diagonal} action
\begin{equation}
\label{eq:accaodiagonal}
g \cdot (v_{1}\otimes \ldots \otimes v_{n})= g \cdot v_{1}\otimes \ldots \otimes g \cdot v_{n},
\end{equation}
for all $v_{1}, \ldots , v_{n}\in V$ and $g \in G_{d}$. 

Hence, $\otimes^{n} V$ is both a left $\k G_{d}$-module and a right $\k S_{n}$-module via place permutations. It is easily seen that these actions commute. In case $\k = \mathbb{C}$, the following theorem goes back to Schur's famous $1927$ article \cite{Schur1927}. As mentioned in the introduction, the result holds over infinite fields of any characteristic.

\vspace{0.2cm}

\begin{theorem} \label{th:SchurWeylduality} Let $V$ be a $d$-dimensional vector space over $\k$ and regard $\otimes^{n} V$ both as the left diagonal $\k G_{d}$-module and the right $\k S_{n}$-module by place permutations. Let $\rho: \k G_{d} \rightarrow {\rm End}_{\, \k}(\otimes^{n} V)$ and $\hat{\rho}: (\k S_{n})^{\circ} \rightarrow {\rm End}_{\, \k}(\otimes^{n} V)$ be the corresponding representations. Then:
\begin{enumerate}
\item [$(a)$] $\rho(\k G_{d})= {\rm End}_{\, \k S_{n}}(\otimes^{n} V)$;
\item[$(b)$] $\hat{\rho}(\k S_{n}))= {\rm End}_{\, \k G_{d}}(\otimes^{n} V)$;
\item [$(c)$] if $d\geq n$, then $\hat{\rho}$ is injective and thus it induces an isomorphism of $\k$-algebras $$(\k S_{n})^{\circ} \cong {\rm End}_{\, \k G_{d} }(\otimes^{n} V).$$
\end{enumerate}
\end{theorem}

\vspace{0.3cm}


\subsection{The classical Schur algebra}\label{subsec2b}

\vspace{0.3cm}

Throughout this paper, we shall adopt Green's viewpoint \cite{Green1} on the polynomial representations of $G_{d}$ and hence work with Schur algebras (see also \cite{Martin}). All of the results presented in this section are valid for infinite fields of arbitrary characteristic.

For all $1\leq i,j \leq d$, let $c_{i,j}: G_{d}\rightarrow \k$ be the previously defined $(i,j)$-th coordinate function. The polynomial algebra $\mathcal{A}(d) \equiv \k[c_{i,j}: 1\leq i,j \leq d]$ has the structure of a bialgebra with comultiplication and counit given by
$$\Delta(c_{i,j}) = \sum_{k=1}^{d} c_{i,k}\otimes c_{k,j} \ {\rm \ and \ } \ \epsilon(c_{i,j})=\delta_{i,j}.$$

For each non-negative integer $n$, let $\mathcal{A}(d,n)$ be the $\k$-subspace of homogeneous polynomials of degree $n$. The bialgebra $\mathcal{A}(d)$ has a natural grading as 
$$\mathcal{A}(d)= \displaystyle{\bigoplus_{n\geq 0} \mathcal{A}(d,n)}.$$
Since each $\mathcal{A}(d,n)$ has the structure of a subcoalgebra of $\mathcal{A}(d)$, its linear dual $\mathcal{A}(d,n)^{\ast}={\rm Hom}_{\, \k}(\mathcal{A}(d,n); \k)$ is an associative algebra over $\k$ which we denote $\mathcal{S}_{\, \scriptscriptstyle{\k}}(d,n)$ and call the (classical) {\it Schur algebra}.

We shall need some notation. If $n\geq 1$, recall that $\mathbf{n} = \{1\,, \ldots \,, n\}$. We write $\Gamma_{\mathbf{n}}(d)$ for the set of all maps $\alpha : \mathbf{n}\rightarrow \mathbf{d}$, identifying each $\alpha$ with the $n$-tuple $(\alpha (1), \ldots , \alpha (n))$ or, equivalently, $(\alpha_{1}, \ldots , \alpha_{n})$. The symmetric group $S_{n}$ acts on the right on $\Gamma_{\mathbf{n}}(d)$ by the rule $\alpha\sigma = (\alpha_{\sigma(1)}, \ldots , \alpha_{\sigma(n)})$, for $\alpha = (\alpha_{1}, \ldots , \alpha_{n})\in \Gamma_{\mathbf{n}}(d)$ and $\sigma \in S_{n}$. Similarly, $S_{n}$ acts on the right on $\Gamma_{\mathbf{n}}(d)\times \Gamma_{\mathbf{n}}(d)$ by $(\alpha, \beta)\sigma = (\alpha\sigma, \beta\sigma)$. If $\alpha, \beta, \gamma, \nu \in \Gamma_{\mathbf{n}}(d)$, $\alpha \sim \beta$ means that $\alpha$ and $\beta$ are in the same $S_{n}$-orbit of $\Gamma_{\mathbf{n}}(d)$ and $(\alpha, \beta)\sim(\gamma, \nu)$ means that $(\alpha, \beta)$ and $(\gamma, \nu)$ are in the same $S_{n}$-orbit of $\Gamma_{\mathbf{n}}(d)\times \Gamma_{\mathbf{n}}(d)$.

The space $\mathcal{A}(d,n)$ has as $\k$-basis the set of all monomials of degree $n$ in the $d^{2}$ variables $c_{i, j}$. Each such monomial can be written as $c_{\alpha, \beta}=c_{\alpha_{1}, \beta_{1}}\ldots c_{\alpha_{n}, \beta_{n}}$ for some $\alpha, \beta \in \Gamma_{\mathbf{n}}(d)$ and there is an ``equality'' rule \cite[Equation~(2.1b)]{Green1}
\begin{equation}
\label{eq:igualdadecalfabeta}
c_{\alpha, \beta}=c_{\gamma, \nu} \ {\rm \ if \ and \ only \ if \ } \ (\alpha, \beta)\sim (\gamma, \nu),
\end{equation}
for all $\alpha, \beta, \gamma, \nu \in \Gamma_{\mathbf{n}}(d)$. Let $\Omega_{n}$ be an arbitrary set of representatives of the $S_{n}$-orbits of $\Gamma_{\mathbf{n}}(d)\times \Gamma_{\mathbf{n}}(d)$. Hence, $\{c_{\alpha, \beta} : (\alpha, \beta) \in \Omega_{n}\}$ is an $\k$-basis of $\mathcal{A}(d,n)$ and thus 
\begin{equation}
\label{eq:dimA(d,r)}
\dim_{\, \scriptscriptstyle{\k}} \left ( \mathcal{A}(d,n) \right )= \binom{d^{2}+n-1}{n}.
\end{equation}
We write $\{\xi_{\alpha, \beta} : (\alpha, \beta) \in \Omega_{n}\}$ for the $\k$-basis of $\mathcal{S}_{\, \scriptscriptstyle{\k}}(d,n)$ dual to that of $\mathcal{A}(d,n)$. Of course, $\dim_{\, \scriptscriptstyle{\k}} ( \mathcal{S}_{\, \scriptscriptstyle{\k}}(d,n) ) = \dim_{\, \scriptscriptstyle{\k}} ( \mathcal{A}(d,n))$ and, for all $\alpha, \beta, \gamma, \nu \in \Gamma_{\mathbf{n}}(d)$, 
\begin{equation}
\label{eq:igualdadecsialfabeta}
\xi_{\alpha, \beta}=\xi_{\gamma, \nu} \ {\rm \ if \ and \ only \ if \ } \ (\alpha, \beta)\sim (\gamma, \nu).
\end{equation}

The algebra structure on $\mathcal{S}_{\, \scriptscriptstyle{\k}}(d,n)$ is the dual of the coalgebra structure on $\mathcal{A}(d,n)$. This implies the following rule for the multiplication in $\mathcal{S}_{\, \scriptscriptstyle{\k}}(d,n)$,
\begin{equation*}
\label{eq:multiplicaremsdr}
\xi\eta (c_{\alpha, \beta})=\sum_{\gamma\in \Gamma_{\mathbf{n}}(d)} \xi(c_{\alpha, \gamma})\eta(c_{\gamma, \beta}),
\end{equation*}
for all $\alpha, \beta \in \Gamma_{\mathbf{n}}(d)$ and $\xi, \eta\in \mathcal{S}_{\, \scriptscriptstyle{\k}}(d,n)$. 

The importance of Schur algebras stems from the fact that they determine the finite-dimensional polynomial representations of $\k G_{d}$. We say that a representation $\rho: \k G_{d} \rightarrow {\rm End}_{\, \k} (W)$ on an $\k$-space $W$ is {\it polynomial} if its  coefficient functions lie in $\mathcal{A}(d)$ and {\it homogeneous of degree} $n$ if its coefficient functions lie in $\mathcal{A}(d,n)$. Following Schur's arguments for $\k=\mathbb{C}$, J. A. Green proved that every polynomial representation of $\k G_{d}$ is a direct sum of homogeneous ones \cite[Theorem $2.2c$]{Green1} and that the category of $\mathcal{S}_{\, \scriptscriptstyle{\k}}(d,n)$-modules is equivalent to the category of homogeneous polynomial modules of $\k G_{d}$ of degree $n$ (see \cite[pp. $23-25$]{Green1}).

For our purposes, we focus on the $\k G_{d}$-module $\otimes^{n} V$ with action given by Expression (\ref{eq:accaodiagonal}). If $\{e_{1}, \ldots , e_{d}\}$ is an $\k$-basis of $V$, then $\{e_{\alpha}^{\otimes} = e_{\alpha(1)}\otimes \ldots \otimes e_{\alpha(n)} : \alpha \in \Gamma_{\mathbf{n}}(d)\}$ is an $\k$-basis of $\otimes^{n} V$. Hence, the left diagonal $\k G_{d}$-action can be expressed as
\begin{equation*}
\label{eq:accaoGld}
g \cdot e_{\beta}^{\otimes} = \sum_{\alpha \in \Gamma_{\mathbf{n}}(d)} c_{\alpha, \beta}(g) e_{\alpha}^{\otimes}, \ {\rm \ for \ } \ g\in G_{d} \ {\rm \ and \ } \ \beta\in \Gamma_{\mathbf{n}}(d).
\end{equation*}

Since the $c_{\alpha, \beta}$ all lie in $\mathcal{A}(d,n)$, this amounts to saying that $\otimes^{n} V$ is a polynomial $\k G_{d}$-module which is homogeneous of degree $n$. As a left $\mathcal{S}_{\, \scriptscriptstyle{\k}}(d,n)$-module,
\begin{equation}
\label{eq:accaoSdrtensores}
\xi e_{\beta}^{\otimes} = \sum_{\alpha \in \Gamma_{\mathbf{n}}(d)} \xi (c_{\alpha, \beta}) e_{\alpha}^{\otimes}, \ {\rm \ for \ } \ \xi \in \mathcal{S}_{\, \scriptscriptstyle{\k }}(d,n) \ {\rm \ and \ } \ \beta \in \Gamma_{\mathbf{n}}(d).
\end{equation}

It is easily seen that the previous $\mathcal{S}_{\, \scriptscriptstyle{\k}}(d,n)$-action commutes with the right action of $S_{n}$ on $\otimes^{n}V$ given by place permutations. In truth, the Schur--Weyl duality between $G_{d}$ and $S_{n}$ on $\otimes^{n} V$ can be stated in terms of $\mathcal{S}_{\, \scriptscriptstyle{\k}}(d,n)$. The proof of the first isomorphism exhibited in the following theorem can be found in \cite[Theorem $(2.6c)$]{Green1}. The other isomorphism follows from \cite[p. $209$, Lemma]{CarterLusztig} and \cite[Section $2.4$]{Green1}.

\vspace{0.2cm}

\begin{theorem}
\label{th:SchurWeyldualitySdr} Let $V$ be a $d$-dimensional vector space over $\k$. Regard $\otimes^{n} V$ both as a left $\mathcal{S}_{\, \scriptscriptstyle{\k}}(d,n)$-module with action given by Equality (\ref{eq:accaoSdrtensores}) and a right $\k S_{n}$-module with action given by Equality (\ref{eq:placepermutsn}). If $d\geq n$, then each action generates the full centralizer of the other on ${\rm End}_{\, \k}(\otimes^{n} V)$ and hence, as $\k$-algebras,
\begin{equation*}
\label{eq:swdschuralgsn}
\mathcal{S}_{\, \scriptscriptstyle{\k}}(d,n) \cong {\rm End}_{\, \k S_{n}}(\otimes^{n} V) \ \  {\rm and} \ \ (\k S_{n})^{\circ} \cong {\rm End}_{\, \mathcal{S}_{\, \scriptscriptstyle{\k}}(d,n)}(\otimes^{n} V).
\end{equation*}
\end{theorem}

\vspace{0.3cm}


\subsection{Representations of the rook monoid}\label{secsec:rookmonoid}

\vspace{0.3cm}

The representation theory of finite inverse semigroups was established in the fifties by W. D. Munn \cite{MunnSemigAlg, MunnMatrixRep, MunnCharact} and I. S. Ponizovski\u{i} \cite{Poniz}. For the special case of the rook monoid, their results were furthered and deepened in zero characteristic by L. Solomon \cite{Solomon}.

Recall that $\mathbf{n}=\{1,\ldots ,n\}$. Our convention for the multiplication in $R_{n}$ is that the composition $\sigma\tau$ of the elements $\sigma, \tau\in R_{n}$ is defined by first applying $\tau$ and then $\sigma$. If $\sigma\in R_{n}$, we write $D(\sigma)\subseteq \mathbf{n}$ for the domain of $\sigma$ and $R(\sigma)\subseteq \mathbf{n}$ for its range.

For instance, $\sigma =
\left ( \begin{array}{ccccc}
1 & 2 & 3 & 4 & 5 \\
- & 4 & 5 & 2 & -
\end{array} \right ) $ and $\tau =
\left ( \begin{array}{ccccc}
1 & 2 & 3 & 4 & 5 \\
3 & 1 & 2 & - & -
\end{array} \right ) $ are elements of $R_{5}$ with $D(\sigma)=\{2,3,4\}$, $R(\sigma)=\{2,4,5\}$ and $D(\tau) = R(\tau)=\{1,2,3\}$.

We agree that $R_{n}$ contains a map $\epsilon_{\emptyset}$ with empty domain and range which behaves as a zero element in $R_{n}$. With this convention, it is easy to see that 
$$|R_{n}|=\sum_{r=0}^{n} \binom{n}{r}^{2} r!$$

If $\sigma \in R_{n}$, we define the {\it rank} of $\sigma$, denoted $rk(\sigma)$, as the size of its domain. We adopt the convention that the only element in $R_{n}$ of rank zero is $\epsilon_{\emptyset}$ and that $S_{0}=\{\epsilon_{\emptyset}\}$ is a group with a single element. Note that any $\sigma \in R_{n}$ of rank $n$ is a permutation of $\mathbf{n}$ and thus $S_{n}\subseteq R_{n}$. A minute's thought reveals that $S_{r}\subseteq R_{n}$ , for $r=0,1, \ldots , n$.

The set of idempotents of $R_{n}$ is the commutative submonoid of all the partial identities $\epsilon_{X}:X\rightarrow X$, with $X\subseteq \mathbf{n}$. For each $X\subseteq \mathbf{n}$, $\epsilon_{X}R_{n}\epsilon_{X}$ is a monoid with identity $\epsilon_{X}$, whose group of units $G_{X}$ is called the {\it maximal subgroup} of $R_{n}$ at $\epsilon_{X}$. 

Two idempotents $\epsilon_{X}$ and $\epsilon_{Y}$ are said to be {\it isomorphic} if $G_{X}\cong G_{Y}$ as groups. If $0\leq r\leq n$ and $X\subseteq \mathbf{n}$ is a set of size $r$, it is not hard to verify that 
$$G_{X}=\{\sigma \in R_{n} : D(\sigma)=R(\sigma)=X\}\cong S_{r}.$$
Thus, each maximal subgroup of $R_{n}$ can be identified with some symmetric group $S_{r}$.

In order to classify the isomorphisms classes of simple modules of $\k R_{n}$, we shall need some notation. If $\sigma\in R_{n}$, we define the {\it inverse} of $\sigma$, denoted $\sigma^{-}$, as the only element of $R_{n}$ which satisfies $D(\sigma^{-})=R(\sigma)$ and $\sigma^{-}\sigma=\epsilon_{D(\sigma)}$.

If $0\leq r\leq n$ and $X\subseteq \mathbf{n}$ is a set of size $r$, we denote by $\iota_{X}$ the unique order-preserving bijection between $\mathbf{r}$ and $X$ (identifying $\mathbf{0}$ with the empty set). Note that any $\sigma \in R_{n}$ of rank $r$ with $D(\sigma)=X$ and $R(\sigma)=Y$ can be mapped to $S_{r}$ via 
$$\mathfrak{p}(\sigma)=\iota_{Y}^{-}\sigma \iota_{X}.$$
In particular, $\mathfrak{p}(\epsilon_{X})=\epsilon_{\mathbf{r}}\in S_{r}$, for all $X\subseteq \mathbf{n}$ of size $r$.

For our purposes, we also need to introduce special algebras. If $1\leq r \leq n$, let $\mathcal{M}_{\binom{n}{r}}(\k S_{r})$ be the $\k$-algebra of all matrices with rows and columns indexed by subsets $I, J$ of $\mathbf{n}$ of size $r$ and entries in $\k S_{r}$. If $r=0$, we identify $\mathcal{M}_{\binom{n}{r}}(\k S_{r})$ with $\k$. Set
\begin{equation}
\label{eq:algebramatrizes}
\mathcal{R}_{n} = \bigoplus_{r=0}^{n}  \mathcal{M}_{\binom{n}{r}}(\k S_{r}).
\end{equation}

If $1\leq r\leq n$ and $I,J\subseteq \mathbf{n}$ are such that $|I|=|J|=r$, let $E_{I,J}$ be the standard matrix with $\epsilon_{\mathbf{r}}\in S_{r}$ in position $(I,J)$ and zeros elsewhere. If $r=0$, set $1_{\k}=E_{\emptyset, \emptyset}$. It is clear that $\mathcal{R}_{n}$ has as $\k$-basis 
\begin{equation}
\label{eq:basealgebramat}
\bigcup_{r=0}^{n} \{\sigma E_{I,J} : \sigma \in S_{r}, I, J\subseteq \mathbf{n}, |I|=|J|=r\}.
\end{equation}

The following result is essentially due to W. D. Munn \cite{MunnSemigAlg,MunnMatrixRep} and I. S. Ponizovski\u{i} \cite{Poniz} although it can be found explicitly in \cite[Lemma $2.17$]{Solomon} and \cite[Theorem $4.4$]{SteinbergMobiusI}.

\vspace{0.2cm}

\begin{theorem}\label{teo:rnisomorfoalg} Let $\mathbb{F}$ be a field of characteristic zero. There is an isomorphism of $\k$-algebras $\phi : \k R_{n}\rightarrow \mathcal{R}_{n}$ given by 
\begin{equation*}\label{eq:imagemfi}
\phi(\sigma)=\sum_{X\subseteq D(\sigma)}\mathfrak{p}(\sigma \epsilon_{X}) E_{\sigma(X), X}
\end{equation*}
with inverse given on basis elements $\sigma \in S_{r}$ and $E_{I,J}$ with $|I|=|J|=r$ by
\begin{equation*}
\label{eq:imagemfimenosum}
\phi^{-1}(\sigma E_{I,J})=\sum_{X\subseteq J} (-1)^{|J|-|X|} (\iota_{I}\sigma \iota_{J}^{-})\epsilon_{X}.
\end{equation*}
In particular, $\k R_{n}$ is a finite-dimensional (split) semisimple algebra over $\k$.
\end{theorem}

\vspace{0.2cm}

As a consequence of the previous theorem, the isomorphisms classes of simple modules of $R_{n}$ are in one-to-one correspondence with those of its maximal subgroups. In order to highlight the constructive aspect of Theorem \ref{teo:rnisomorfoalg}, we express this result in terms of matrix representations of $\k R_{n}$.

Let $\mu\vdash r$ with $0\leq r\leq n$ and let $\rho_{\mu}: \k S_{r} \rightarrow \mathcal{M}_{k}(\k)$ be an irreducible matrix representation of $\k S_{r}$. If $r=0$, we agree that there is an empty partition $\mu=(0)$ and the corresponding irreducible representation of $\k S_{0} $ is given by $\rho_{(0)}(\epsilon_{\emptyset})=1_{\k}$.

It follows from Theorem \ref{teo:rnisomorfoalg} that $\rho_{\mu}$ induces a matrix representation of $\k R_{n}$, denoted by $\rho_{\mu}^{\ast}: \k R_{n} \rightarrow \mathcal{M}_{k\,\binom{n}{r}}(\k)$, and given by
\begin{equation}
\label{eq:repirredurhoestrela}
\rho_{\mu}^{\ast}(\sigma) = \sum_{\begin{smallmatrix} X\subseteq \mathbf{n}, |X|=r, \\ rk(\sigma\epsilon_{X})=r \end{smallmatrix}} \rho_{\mu}(\mathfrak{p}(\sigma\epsilon_{X})) E_{R(\sigma\epsilon_{X}), D(\sigma\epsilon_{X})},
\end{equation}
for all $\sigma \in R_{n}$. The previous expression is due to L. Solomon (see Equation $2.26$ in \cite{Solomon}). In the setting of finite inverse semigroups, W. D. Munn \cite{MunnMatrixRep} showed that these representations determine the isomorphism type of $\k R_{n}$.

\vspace{0.2cm}

\begin{theorem} [Munn]
\label{teo:irredutiveisRn}  Let $\k$ be a field of characteristic zero. The set 
$$\{\rho_{\mu}^{\ast} : \mu\vdash r, \ {\rm for \ all \ } 0\leq r\leq n\}$$
is a complete set of inequivalent irreducible matrix representations of $\k R_{n}$. Thus, the isomorphism classes of simple modules of $\k R_{n}$ are indexed by the set of partitions
 $$\{\mu\vdash r : \ 0\leq r\leq n\}$$
\end{theorem}

\vspace{0.2cm}


\vspace{0.3cm}

\section{Representing the rook monoid and the extended Schur algebra on tensors}\label{sec2}

\vspace{0.3cm}

Throughout this section, $\k$ is a field of characteristic zero, $d$ and $n$ are positive integers such that $d\geq n$ and $G_{d}$ is viewed as the subgroup of $G_{d+1}$ of all matrices of the form 
$$\left [
\begin{array}{cc}
g & 0  \\
0 & 1
\end{array} \right ],$$ 
where $g \in G_{d}$.
All the results exposed in sections \ref{secsec:restriction} and \ref{secsec:ExtendedSchurAlg} remain valid if $\k$ is replaced by an infinite field of any characteristic.

\vspace{0.3cm}

\subsection{A natural restriction}\label{secsec:restriction}

\vspace{0.3cm}

Let $U=V\oplus U_{0}$ be an $\k$-space of dimension $d+1$, where $V$ and $U_{0}$ are seen as subspaces of $U$ of respective dimensions $d$ and $1$. As before, $\{e_{1}, \ldots , e_{d}\}$ is an arbitrary but fixed basis of $V$ and $U_{0} = \k e_{\infty}$ for some vector $e_{\infty}\in U$ which turns $\{e_{1}, \ldots , e_{d}, e_{\infty}\}$ into an $\k$-basis of $U$. We assume the linear ordering $1<2<\ldots <d<\infty$.

It follows from Section \ref{sec:prelim} that $\otimes^{n} U$ is a finite-dimensional homogeneous polynomial $\k G_{d+1}$-module of degree $n$ via the diagonal action of $G_{d+1}$ on $\otimes^{n} U$ (see Equation (\ref{eq:accaodiagonal}) with $d$ replaced by $d+1$). We study this action's restriction to $G_{d}$.
 
To do so, we need some notation. Let $1\leq r \leq n$ and let $X=\{x_{1}< \ldots < x_{r}\}\subseteq \mathbf{n}$ be a set of size $r$. We write $\Gamma_{X}(d)$ for the set of all maps $\alpha : X \rightarrow \mathbf{d}$ identifying $\alpha \in \Gamma_{X}(d)$ with the $r$-tuple $(\alpha(x_{1}), \ldots , \alpha(x_{r}))\in \mathbf{d}^{r}$. We also agree that $\Gamma_{\emptyset}(d)$ has a single element which is identified with the empty set.

If $\alpha, \beta \in \Gamma_{X}(d)$, then $c_{\alpha, \beta}=c_{\alpha(x_{1}), \beta(x_{1})}\ldots c_{\alpha(x_{r}), \beta(x_{r})}\in \mathcal{A}(d,r)$ is as before a monomial of degree $r$ in the $d^{2}$ variables $c_{i,j}$. If $\iota_{X}$ is the unique order-preserving bijection between $\mathbf{r}$ and X (see Section \ref{secsec:rookmonoid}), the same monomial can be written as
\begin{equation}
\label{eq:translacaocalfabeta}
c_{\alpha, \beta}=c_{\alpha \iota_{X}, \beta \iota_{X}}= c_{\alpha(\iota_{X}(1)), \beta(\iota_{X}(1))}\ldots c_{\alpha (\iota_{X}(r)), \beta(\iota_{X}(r))} = c_{\gamma, \nu},
\end{equation}
where $\gamma=\alpha \iota_{X}, \nu=\beta \iota_{X}\in \Gamma_{\mathbf{r}}(d)$.

This notation is particularly useful to represent the $\k$-basis of $\otimes^{n} U$ induced by $\{e_{1}, \ldots , e_{d}, e_{\infty}\}$. If $X\subseteq \mathbf{n}$ and $\alpha \in \Gamma_{X}(d)$, the decomposable tensor $e_{\alpha}^{\otimes}\in \otimes^{n} U$ is defined as $e_{\alpha}^{\otimes} = e_{\widehat{\alpha}(1)} \otimes \ldots \otimes e_{\widehat{\alpha}(n)}$, where the map $\widehat{\alpha} : \mathbf{n} \rightarrow \mathbf{d}\cup \{\infty\}$ is given by
\begin{displaymath}
\widehat{\alpha} (i) =
 \begin{cases}
\alpha(i) & \textrm{if $i \in X$,} \\
\infty & \textrm{otherwise.}
 \end{cases}\
\end{displaymath}

It is clear that $\otimes^{n} U$ has as $\k$-basis the set $\{e_{\alpha}^{\otimes} : \alpha \in \Gamma_{X}(d), \ X\subseteq \mathbf{n} \}$. For instance, if $d=6$, $n=5$, $X=\{1,3,5\}$ and $\alpha=(6,2,2)\in \Gamma_{\{1,3,5\}}(6)$, then the corresponding basis element of $\otimes^{5} U$ is given by
$$e_{\alpha}^{\otimes} = e_{6}\otimes e_{\infty} \otimes e_{2} \otimes e_{\infty} \otimes e_{2} \in \otimes^{5} U.$$

\vspace{0.3cm}

\begin{proposition} 
\label{prop:accaoreduzida} 
Let $U$ be a $(d+1)$-dimensional vector space over $\k$. The restriction to $G_{d}$ of the left diagonal action of $G_{d+1}$ on $\otimes^{n} U$ is given by
\begin{equation}
\label{eq:accaoGreduzida}
g \cdot e_{\beta}^{\otimes} = \sum_{\alpha \in \Gamma_{X}(d)} c_{\alpha , \beta}(g) e_{\alpha}^{\otimes}, {\rm \ for \ all } \ g\in G_{d} , \beta \in \Gamma_{X}(d) {\rm \ and \ } X\subseteq \mathbf{n}.
\end{equation}
\end{proposition}
\begin{proof} Let $g\in G_{d}$. Under the identification $G_{d} \leq G_{d+1}$, the natural action of $g$ on basis elements of $U$ becomes
$$g \cdot e_{j} = \sum_{i=1}^{d} c_{i,j}(g) e_{i} + 0 \, . e_{\infty}\,, \ {\rm for} \ j=1,\ldots, d\,, \ {\rm and} \  g \cdot e_{\infty} = e_{\infty}.$$

We turn to the action of $g$ on basis elements of the diagonal $\k G_{d+1}$-module $\otimes^{n} U$. If $X=\emptyset$, then we have $g \cdot e^{\otimes}_{\emptyset}= g \cdot e_{\infty} \otimes \ldots \otimes g \cdot e_{\infty} = e^{\otimes}_{\emptyset}$. Let $1\leq r\leq n$ and let $X=\{x_{1}< \ldots < x_{r}\}\subseteq \mathbf{n}$ be a set of size $r$. If $\beta\in \Gamma_{X}(d)$, the corresponding basis element of $\otimes^{n} U$ is $e_{\beta}^{\otimes}= e_{\widehat{\beta}(1)}\otimes \ldots \otimes e_{\widehat{\beta}(n)}$, where $\widehat{\beta}: \mathbf{n} \rightarrow \mathbf{d}\cup\{\infty\}$ is defined as above. It follows that $g \cdot e_{\widehat{\beta}(i)} = e_{\infty}$, for all $i\notin X$ and also that \linebreak $g \cdot e_{\beta}^{\otimes} = \displaystyle{g \cdot e_{\widehat{\beta}(1)} \otimes \ldots \otimes g \cdot e_{\widehat{\beta}(n)} \in \sum_{\alpha\in \Gamma_{X}(d)} a_{\alpha} e_{\alpha}^{\otimes}}$, for some $a_{\alpha}\in \k$. It is now easy to see that 
$$g \cdot e_{\beta}^{\otimes} = \sum_{\alpha\in \Gamma_{X}(d)}(c_{\alpha(x_{1}), \beta(x_{1})}(g)\ldots c_{\alpha(x_{r}),\beta(x_{r})}(g)) e_{\alpha}^{\otimes}= \sum_{\alpha\in\Gamma_{X}(d)} c_{\alpha, \beta}(g) e_{\alpha}^{\otimes},$$ where $c_{\alpha, \beta}\in \mathcal{A}(d,r)$, for all $\alpha \in \Gamma_{X}(d)$.
\end{proof}

\vspace{0.3cm}


\subsection{The extended Schur algebra}\label{secsec:ExtendedSchurAlg}

\vspace{0.3cm}

The coefficient space produced by the action of $G_{d}$ on $\otimes^{n} U$ described in Proposition \ref{prop:accaoreduzida} suggests that $\otimes^{n} U$ can be seen as a representation of a special Schur algebra.

Let $\mathcal{A}_{\mathbf{n}}(d) \, = \, < c_{\alpha , \beta} : (\alpha, \beta) \in \Gamma_{X}(d)\times \Gamma_{X}(d), X\subseteq \mathbf{n} >$ be the $\k$-space spanned by all the monomials $c_{\alpha , \beta}$ of degree at most $n$ in the variables $c_{i,j}$. A moment's thought and Equation (\ref{eq:translacaocalfabeta}) reveal that $\mathcal{A}_{\mathbf{n}}(d)$ is the direct sum of the first $n+1$ homogeneous $\k$-spaces $\mathcal{A}(d,r)$ of the graded bialgebra $\mathcal{A}(d)$. 

It follows that $\mathcal{A}_{\mathbf{n}}(d)$ has as $\k$-basis the set of all distinct monomials of degree $r$ with $0\leq r\leq n$. By Equation (\ref{eq:igualdadecalfabeta}), we index this basis with the set $ \Omega_{0} \cup \Omega_{1} \cup \ldots \cup \Omega_{n}$, where $\Omega_{0}=\Gamma_{\emptyset}(d)\times \Gamma_{\emptyset}(d)$ and each $\Omega_{r}$, with $1\leq r \leq n$, is a set of representatives of the $S_{r}$-orbits of $\Gamma_{\mathbf{r}}(d)\times \Gamma_{\mathbf{r}}(d)$.

As mentioned previously, each $\mathcal{A}(d,r)$ is a subcoalgebra of $\mathcal{A}(d)$. Thus, $\mathcal{A}_{\mathbf{n}}(d)$ enherits a coalgebra structure with comultiplication and counit given, respectively, by
\begin{equation}
\label{eq:coalgebraA}
\Delta(c_{\alpha, \beta}) = \sum_{\gamma\in \Gamma_{X}(d)} c_{\alpha , \gamma} \otimes c_{\gamma , \beta} \ {\rm \ and \ } \  \epsilon(c_{\alpha , \beta}) = \delta_{\alpha , \beta},
\end{equation}
for all $\alpha, \beta \in \Gamma_{X}(d)$ and $X\subseteq \mathbf{n}$. By a standard fact, it follows that 
$$\mathcal{A}_{\mathbf{n}}(d)^{\ast}={\rm Hom}_{\k} \left (\mathcal{A}_{\mathbf{n}}(d); \k \right )$$
is an associative $\k$-algebra of finite dimension.

\vspace{0.2cm}

\begin{definition}
\label{def:extendedSchurAlg} Let $d$ and $n$ be positive integers. The extended Schur algebra for $d$ and $n$ over the field $\k$, denoted $\mathcal{S}_{\, \scriptscriptstyle{\k}}(d, \mathbf{n})$, is the associative $\k$-algebra of finite dimension given by the linear dual $\mathcal{S}_{\, \scriptscriptstyle{\k}}(d, \mathbf{n}) = \mathcal{A}_{\mathbf{n}}(d)^{\ast}={\rm Hom}_{\k} \left (\mathcal{A}_{\mathbf{n}}(d); \k \right )$.
\end{definition}

\vspace{0.2cm}

Since $\mathcal{A}_{\mathbf{n}}(d) = \displaystyle{\bigoplus_{r=0}^{n} \mathcal{A}(d,r)}$, it follows that $\mathcal{S}_{\, \scriptscriptstyle{\k}}(d, \mathbf{n})$ can be regarded as the $\k$-algebra 
$$\displaystyle{\bigoplus_{r=0}^{n}  \mathcal{S}_{\, \scriptscriptstyle{\k}}(d,r)},$$
where each $\mathcal{S}_{\, \scriptscriptstyle{\k}}(d,r) $ is a classical Schur algebra. Indeed, if $\xi \in  \mathcal{S}_{\, \scriptscriptstyle{\k}}(d,r) $, we identify $\xi$ with an element of $\mathcal{S}_{\, \scriptscriptstyle{\k}}(d, \mathbf{n}) = \displaystyle{\bigoplus_{r=0}^{n} \mathcal{A}(d,r)^{\ast}}$ by making it zero on all monomials whose degree is different from $r$. Under this identification, 
$$\displaystyle{\bigcup_{r=0}^{n}\{\xi_{\alpha , \beta} : (\alpha, \beta) \in \Omega_{r}\}}$$
is the $\k$-basis of $\mathcal{S}_{\, \scriptscriptstyle{\k}}(d, \mathbf{n})$ which is dual to the $\k$-basis of $\mathcal{A}_{\mathbf{n}}(d)$ given by 
$$\displaystyle{\bigcup_{r=0}^{n}\{c_{\alpha , \beta} : (\alpha, \beta) \in \Omega_{r}\}}.$$
Combined with Equation (\ref{eq:dimA(d,r)}), this implies that $\dim_{\, \scriptscriptstyle{\k}} \left ( \mathcal{S}_{\, \scriptscriptstyle{\k}}(d, \mathbf{n}) \right ) = \displaystyle{\binom{d^{2}+n}{n}}$. Indeed,
$$ \sum_{r=0}^{n} \dim_{\, \scriptscriptstyle{\k}} \left ( \mathcal{S}_{\, \scriptscriptstyle{\k}}(d,r) \right ) = \sum_{r=0}^{n} \dim_{\, \scriptscriptstyle{\k}} \left ( \mathcal{A}(d,r) \right ) = \sum_{r=0}^{n} \binom{d^{2}+r-1}{r} = \binom{d^{2}+n}{n}\,.$$
If $0\leq r\leq n$ and $\alpha, \beta \in \Gamma_{\mathbf{r}}(d)$, $\xi_{\alpha , \beta}$ is the element of $\mathcal{S}_{\, \scriptscriptstyle{\k}}(d, \mathbf{n})$ given by
\begin{equation}
\label{eq:basealgebraS}
\xi_{\alpha , \beta}(c_{\gamma, \nu}) 
 = 
 \begin{cases}
1 & \textrm{if $r=k$ and $(\alpha, \beta)\sim (\gamma, \nu)$} \\
0 & \textrm{otherwise,}
 \end{cases}
\end{equation}
for all $\gamma, \nu \in \Gamma_{\mathbf{k}}(d)$ and $0\leq k\leq n$. As with the classical case (see Equation (\ref{eq:igualdadecsialfabeta})), if $0\leq k\leq n$ and $\kappa, \varsigma \in \Gamma_{\mathbf{k}}(d)$, we also have an equality rule to take into account, namely,
\begin{equation*}
\label{eq:igualdadeemS}
\xi_{\alpha , \beta}= \xi_{\kappa , \varsigma} \ {\rm \ if \ and \ only \ if \ } \ r=k \ {\rm \ and \ } \ (\alpha, \beta)\sim (\kappa, \varsigma).
\end{equation*}

For now, it is enough to work with the index set $\displaystyle{\bigcup_{r=0}^{n} \Gamma_{\mathbf{r}}(d)}$, since it follows from equations (\ref{eq:translacaocalfabeta}) and (\ref{eq:basealgebraS}) that, for $\alpha, \beta \in \Gamma_{X}(d)$ with $X\subseteq \mathbf{n}$ such that $|X|=r$,
\begin{equation}
\label{eq:translacaocsialfabeta}
\xi_{\alpha, \beta}= \xi_{\alpha \iota_{X}, \beta \iota_{X}} = \xi_{\gamma, \nu} \in \mathcal{S}_{\,  \scriptscriptstyle{\k}}(d, \mathbf{n}),
\end{equation}
with $\gamma=\alpha \iota_{X}, \nu=\beta \iota_{X} \in \Gamma_{\mathbf{r}}(d)$ and $\iota_{X}:\mathbf{r}\rightarrow X$ defined as before.

The multiplication in $\mathcal{S}_{\, \scriptscriptstyle{\k}}(d, \mathbf{n})$ follows from the coalgebra structure on $\mathcal{A}_{\mathbf{n}}(d)$ and hence from Equation (\ref{eq:coalgebraA}). Thus, if $\xi , \eta \in \mathcal{S}_{\, \scriptscriptstyle{\k}}(d, \mathbf{n})$, the product $\xi\eta$ is defined on any monomial $c_{\alpha , \beta}$, with $\alpha, \beta \in \Gamma_{X}(d)$ and $X\subseteq \mathbf{n}$, by
\begin{equation}
\label{eq:produtoemS}
\xi\eta (c_{\alpha , \beta})=\sum_{\gamma\in \Gamma_{X}(d)} \xi(c_{\alpha ,\gamma})\eta(c_{\gamma ,\beta}).
\end{equation}

The unit element $\epsilon \in \mathcal{S}_{\, \scriptscriptstyle{\k}}(d, \mathbf{n})$ is given by $\epsilon(c)=c(I_{d})$, for all $c \in  \mathcal{A}_{\mathbf{n}}(d)$. Note that $\epsilon \in \mathcal{S}_{\, \scriptscriptstyle{\k}}(d, \mathbf{n})$ can be expressed as $\epsilon = \displaystyle{\sum_{r=0}^{n} \epsilon_{r}}$, where $\epsilon_{r}$ is the identity of $\mathcal{S}_{\, \scriptscriptstyle{\k}}(d,r)$.

We now turn to the tensor space $\otimes^{n} U$ and show that it can be given the structure of a left $ \mathcal{S}_{\, \scriptscriptstyle{\k}}(d, \mathbf{n})$-module. We have however a stronger statement.

\vspace{0.2cm}

\begin{proposition}
\label{prop:categoriasequiv}
The category of finite-dimensional $\k G_{d}$-modules whose coefficient functions lie in $\displaystyle{\mathcal{A}_{\mathbf{n}}(d)= \bigoplus_{r=0}^{n} \mathcal{A}(d,r)}$ is equivalent to that of $\mathcal{S}_{\, \scriptscriptstyle{\k}}(d, \mathbf{n})$-modules.
\end{proposition}
\begin{proof}
Let $g\in G_{d}$. For any $X\subseteq \mathbf{n}$ and $\alpha, \beta\in \Gamma_{X}(d)$, define $e_{g}(c_{\alpha, \beta})=c_{\alpha, \beta}(g)$. By linear extension, the map $e_{g}: c\mapsto c(g)$ is a well-defined linear homomorphism of $\mathcal{S}_{\, \scriptscriptstyle{\k}}(d, \mathbf{n}) = {\rm Hom}_{\k} \left (\mathcal{A}_{\mathbf{n}}(d); \k \right )$. If $g, g' \in G_{d}$, it is clear from Equation (\ref{eq:produtoemS}) that $e_{g}e_{g '}=e_{g g '}$ and $e_{I_{d}}=\epsilon \in \mathcal{S}_{\, \scriptscriptstyle{\k}}(d, \mathbf{n})$. Hence, the map $e: g \mapsto e_{g}$ can be linearly extended to an $\k$-algebra homomorphism $e: \k G_{d} \rightarrow \mathcal{S}_{\, \scriptscriptstyle{\k}}(d, \mathbf{n})$.

We assume that any map $f: G_{d} \rightarrow \k$ is identified with its unique linear extension $f: \k G_{d} \rightarrow \k$. Under this assumption, if $k\in \k G_{d}$, then $e_{k}: \mathcal{A}_{\mathbf{n}}(d)\rightarrow \k$ is given by $e_{k}(c)=c(k)$, for all $c\in \mathcal{A}_{\mathbf{n}}(d)$. The arguments in \cite[Proposition $2.4b, (i)$]{Green1} apply {\it mutatis mutandis} to $e: \k G_{d} \rightarrow \mathcal{A}_{\mathbf{n}}(d)^{\ast}$ and hence $e$ is surjective.

We now show that $f: \k G_{d} \rightarrow \k$ belongs to $\mathcal{A}_{\mathbf{n}}(d)$ if and only if $f(k)=0$, for all $k\in \ker e$. If $f\in \mathcal{A}_{\mathbf{n}}(d)$ and $k\in \ker e$, then $e_{k}=0$ and $e_{k}(f)=f(k)=0$. Conversely, let $f: \k G_{d} \rightarrow \k$ be such that $f(k)=0$, for all $k\in \ker e$. Since $e$ is surjective, for all $\xi\in \mathcal{S}_{\, \scriptscriptstyle{\k}}(d, \mathbf{n})$, there is some $k \in \k G_{d}$ such that $\xi = e_{k}$. Hence, we define $y\in \mathcal{S}_{\, \scriptscriptstyle{\k}}(d, \mathbf{n})^{\ast}$ by $y(\xi)=y(e_{k})=f(k)$, for all $k\in \k G_{d}$. The condition that $f(k)=0$, for all $k\in \ker e$ ensures that $y$ is a well-defined element of $\mathcal{S}_{\, \scriptscriptstyle{\k}}(d, \mathbf{n})^{\ast}$. Since $\mathcal{A}_{\mathbf{n}}(d)$ is finite-dimensional, we have that $\mathcal{A}_{\mathbf{n}}(d) \cong \mathcal{A}_{\mathbf{n}}(d)^{\ast\ast} = \mathcal{S}_{\,  \scriptscriptstyle{\k}}(d, \mathbf{n})^{\ast}$ and hence there is some $c\in  \mathcal{A}_{\mathbf{n}}(d)$ such that $y=c$. Let $k\in \k G_{d}$, then $f(k)=y(e_{k})=e_{k}(c)=c(k)$ and thus $c=f\in  \mathcal{A}_{\mathbf{n}}(d)$.

Let $V$ be a finite-dimensional left $\k G_{d}$-module with basis $\{v_{b} : b\in B\}$ and associated action $g \cdot  v_{b} = \displaystyle{ \sum_{a\in B} \alpha_{a,b}(g) v_{a}}$, for all $b\in B$ and $g\in G_{d}$. Suppose that $\alpha_{a,b}\in \mathcal{A}_{\mathbf{n}}(d)$, for all $a,b \in B$. Then $\alpha_{a,b}(k)=0$, for all $k\in \ker e$ and all $a,b \in B$. The action 
$$e_{g} \cdot v_{b} = \sum_{a\in B} e_{g}(\alpha_{a,b}) v_{a},$$
for all $g\in G_{d}$, $b\in B$, turns $V$ into a left $\mathcal{S}_{\, \scriptscriptstyle{\k}}(d, \mathbf{n})$-module. Indeed, for all $\xi \in \mathcal{S}_{\, \scriptscriptstyle{\k}}(d, \mathbf{n})$, there is $k\in \k G_{d}$ such that $\xi = e_{k}$. If $\xi = e_{k} = e_{k '}$ for $k ' \in \k G_{d}$, then $k-k ' \in \ker e$ and $\alpha_{a,b}(k-k ')=0$ for all $a,b \in B$. Hence, $e_{k}(\alpha_{a,b})=e_{k '}(\alpha_{a,b})$ for all $a,b\in B$ and the $\mathcal{S}_{\, \scriptscriptstyle{\k}}(d, \mathbf{n})$-action is well defined.

Conversely, if $V$ is a left $\mathcal{S}_{\, \scriptscriptstyle{\k}}(d, \mathbf{n})$-module with basis $\{v_{b} : b\in B\}$ and associated action $\xi \cdot v_{b}= \displaystyle{\sum_{a\in B} \xi(\alpha_{a,b}) v_{a}}$, for all $b\in B$ and all $\xi \in \mathcal{S}_{\, \scriptscriptstyle{\k}}(d, \mathbf{n})$, then $V$ can be viewed as an $\k G_{d}$-module with action given by 
$$g \cdot v_{b} = e_{g} \cdot v_{b} = \sum_{a\in B} \alpha_{a,b}(g),$$
for all $g\in G_{d}$ and all $b\in B$, where $e_{g}=\xi\in \mathcal{S}_{\, \scriptscriptstyle{\k}}(d, \mathbf{n})$ and $\alpha_{a,b}(g)=\xi(\alpha_{a,b})$. Once again, the previous properties show that $\alpha_{a,b}\in  \mathcal{A}_{\mathbf{n}}(d)$ and this action is well-defined. The proof is complete and we can now identify both categories by the simple rule: $k \cdot v = e_{k} \cdot v$, for all $k\in \k G_{d}$ and $v\in V$, where $V$ is an object of either categories.
\end{proof}

\vspace{0.2cm}

It follows from the proof of this result that the left action of $\mathcal{S}_{\, \scriptscriptstyle{\k}}(d, \mathbf{n})$ on $\otimes^{n} U$ is given, for $\xi \in \mathcal{S}_{\, \scriptscriptstyle{\k}}(d, \mathbf{n})$, $\beta \in \Gamma_{X}(d)$ and $X\subseteq \mathbf{n}$, by
\begin{equation} 
\label{eq:accaotensorS}
\xi e_{\beta}^{\otimes}= \sum_{\alpha\in \Gamma_{X}(d)} \xi(c_{\alpha , \beta})e_{\alpha}^{\otimes}.
\end{equation}
 
We end this section with an important fact which follows from the semisimplicity of the classical Schur algebras.

\vspace{0.2cm}

\begin{proposition}
\label{prop:Ssemisimples}
The extended Schur algebra $\mathcal{S}_{\, \scriptscriptstyle{\k}}(d, \mathbf{n})$ is a semisimple algebra over $\k$.
\end{proposition}
\begin{proof} As referred previously, we may regard $\mathcal{S}_{\, \scriptscriptstyle{\k}}(d, \mathbf{n})$ as $\displaystyle{\bigoplus_{r=0}^{n}  \mathcal{S}_{\, \scriptscriptstyle{\k}}(d,r)}$, where $ \mathcal{S}_{\, \scriptscriptstyle{\k}}(d,r)$ is the classical Schur algebra. A proof of the semisimplicity of $ \mathcal{S}_{\, \scriptscriptstyle{\k}}(d,r)$ can be found in \cite[Corollary $(2.6e)$]{Green1} and hence the semisimplicity of $\mathcal{S}_{\, \scriptscriptstyle{\k}}(d, \mathbf{n})$ follows.
\end{proof}

\vspace{0.3cm}

\subsection{Schur--Weyl duality between the rook monoid and the extended Schur algebra}\label{secsec:SchurWeylRn}

\vspace{0.3cm}

In what follows, we describe the centralizer algebra ${\rm End}_{\, \mathcal{S}_{\, \scriptscriptstyle{\k}}(d, \mathbf{n})} (\otimes^{n} U)$ of $\mathcal{S}_{\, \scriptscriptstyle{\k}}(d, \mathbf{n})$ on the left module $\otimes^{n} U$ on which $\mathcal{S}_{\, \scriptscriptstyle{\k}}(d, \mathbf{n})$ acts according to Equation (\ref{eq:accaotensorS}). Throughout this section, $\k$ is a field of characteristic zero.

For each $X\subseteq \mathbf{n}$, we denote by $W_{X}$ the $\k$-subspace of $\otimes^{n} U$ spanned by all the decomposable tensors $e_{\alpha}^{\otimes}$, with $\alpha\in \Gamma_{X}(d)$. Since $\{e_{\alpha}^{\otimes} : \alpha\in \Gamma_{X}(d), X\subseteq \mathbf{n}\}$ is an $\k$-basis of $\otimes^{n} U$, we have the following direct sum decomposition
\begin{equation}
\label{eq:decompWx}
\otimes^{n} U = \bigoplus_{X\subseteq \mathbf{n}} W_{X}.
\end{equation}

It follows from Equation (\ref{eq:accaotensorS}) that the left action of an arbitrary $\xi\in\mathcal{S}_{\, \scriptscriptstyle{\k}}(d, \mathbf{n})$ on $e_{\alpha}^{\otimes}$, with $\alpha\in\Gamma_{X}(d)$, is such that $\xi e_{\alpha}^{\otimes}\in W_{X}$. Hence, $W_{X}$ is an $\mathcal{S}_{\, \scriptscriptstyle{\k}}(d, \mathbf{n})$-submodule of $\otimes^{n} U$ and (\ref{eq:decompWx}) is a decomposition of $\otimes^{n} U$ as a direct sum of left $\mathcal{S}_{\, \scriptscriptstyle{\k}}(d, \mathbf{n})$-submodules.

The next result shows how ${\rm End}_{\, \mathcal{S}_{\, \scriptscriptstyle{\k}}(d, \mathbf{n})} (\otimes^{n} U)$ decomposes into its building blocks.

\vspace{0.2cm}
  \begin{lemma}     \label{lem:distributividadeSmod} Let $U = V\oplus \k e_{\infty}$ be an $\k$-space of dimension $d+1$ such that $V$ is a $d$-dimensional $\k$-subspace of $U$.
\begin{enumerate}
\item [$(a)$] If $0\leq r\leq n$, the tensor space $\otimes^{r} V$ is a left $\mathcal{S}_{\, \scriptstyle{\k}}(d, \mathbf{n})$-module for which there is an isomorphism of $\k$-algebras 
$${\rm End}_{ \mathcal{S}_{\, \scriptstyle{\k}}(d, \mathbf{n})}(\otimes^{r} V) \cong (\k S_{r})^{\circ}.$$
\item [$(b)$] If $0\leq r\leq n$ and $X\subseteq \mathbf{n}$ is a set of size $r$, then, as $\mathcal{S}_{\, \scriptstyle{\k}}(d, \mathbf{n})$-modules,
$W_{X}\cong \otimes^{r} V$.
\item [$(c)$] There is a left $\mathcal{S}_{\, \scriptstyle{\k}}(d, \mathbf{n})$-module isomorphism such that
$$\otimes^{n} U \cong \bigoplus_{r=0}^{n} \binom{n}{r} \otimes^{r} V$$
where $\displaystyle{ \binom{n}{r} \otimes^{r} V}$ means a direct sum of $\binom{n}{r}$ copies of $\otimes^{r} V$.
\end{enumerate}
\end{lemma}
\begin{proof}
$(a)$ If $r=0$, we agree that $\otimes^{0} V=\k$ with trivial left $\mathcal{S}_{\, \scriptscriptstyle{\k}}(d, \mathbf{n})$-action and right $\k S_{0}$-action. If $r\geq 1$, we view $\otimes^{r} V$ both as a left $\mathcal{S}_{\, \scriptscriptstyle{\k}}(d,r)$-module (via Equation (\ref{eq:accaoSdrtensores})) and a right $\k S_{r}$-module with respect to place permutations. The $\mathcal{S}_{\, \scriptscriptstyle{\k}}(d,r)$-action on $\otimes^{r} V$ is easily extended to an $\mathcal{S}_{\, \scriptscriptstyle{\k}}(d, \mathbf{n})$-action by defining $\xi w=0$, for all $w\in \otimes^{r} V$ and all $\xi \in \mathcal{S}_{\, \scriptscriptstyle{\k}}(d,k)$ with $k\neq r$. It follows that ${\rm End}_{\, \mathcal{S}_{\, \scriptscriptstyle{\k}}(d, \mathbf{n})}(\otimes^{r} V) \equiv  {\rm End}_{\, \mathcal{S}_{\, \scriptscriptstyle{\k}}(d,r)}(\otimes^{r} V)$ and, by Theorem \ref{th:SchurWeyldualitySdr}, ${\rm End}_{\, \mathcal{S}_{\, \scriptscriptstyle{\k}}(d, \mathbf{n})}(\otimes^{r} V) \cong (\k S_{r})^{\circ}$. $(b)$ The case $r=0$ is trivial since $W_{\emptyset}=\k e_{\emptyset}^{\otimes}$. Let $r\geq 1$. If $X=\{x_{1} < \ldots < x_{r}\}\subseteq \mathbf{n}$ is of size $r$, the map $T_{X}: W_{X} \rightarrow \otimes^{r} V$ defined by $T_{X}(e_{\alpha}^{\otimes})=e_{\alpha(x_{1})}\otimes \ldots \otimes e_{\alpha(x_{r})}$, for all $\alpha \in \Gamma_{X}(d)$, and extended linearly to $W_{X}$, is easily seen to be a left $\mathcal{S}_{\, \scriptscriptstyle{\k}}(d, \mathbf{n})$-module isomorphism. $(c)$ This follows from $(a)$ and $(b)$ and the direct sum decomposition described in Equality (\ref{eq:decompWx}).
\end{proof}

\vspace{0.2cm}

It is worth noting that the previous lemma remains valid for arbitrary infinite fields. 

Recall that the algebra of matrices 
$$\mathcal{R}_{n} = \bigoplus_{r=0}^{n}  \mathcal{M}_{\binom{n}{r}}(\k S_{r})$$
has an $\k$-basis given by Equality (\ref{eq:basealgebramat}) and that $\mathcal{R}_{n}$ is isomorphic to $\k R_{n}$ (Theorem \ref{teo:rnisomorfoalg}).

\vspace{0.2cm}

\begin{lemma}
\label{le:accaoalgmatrizes} 
Let $U$ be an $\k$-space of dimension $d+1$. Let $\alpha\in \Gamma_{X}(d)$ for some $X\subseteq \mathbf{n}$ and let $e_{\alpha}^{\otimes}$ be the corresponding basis element of $\otimes^{n} U$. If $\sigma \in S_{r}$ and $I, J \subseteq \mathbf{n}$ are such that $|I|=|J|=r$, define
\begin{equation}
\label{eq:accaoalgmatrizes} 
e_{\alpha}^{\otimes} \centerdot (\sigma E_{I, J})  = \delta_{X,I} e^{\otimes}_{\alpha \iota_{I}\sigma \iota_{J}^{-}}, \ {\rm if} \ 1\leq r\leq n \ \ {\rm and} \ \ e_{\alpha}^{\otimes} \centerdot (\sigma E_{I, J}) = \delta_{X, \emptyset}e_{\emptyset}^{\otimes}, \ {\rm if} \ r=0.
\end{equation}
Then Equation (\ref{eq:accaoalgmatrizes}) gives $\otimes^{n} U$ a right $\mathcal{R}_{n}$-module structure for which the action of $\mathcal{R}_{n}$ commutes with that of $\mathcal{S}_{\, \scriptscriptstyle{\k}}(d, \mathbf{n})$ on $\otimes^{n} U$.
\end{lemma}
\begin{proof} Let $X\subseteq \mathbf{n}$ and $\alpha\in \Gamma_{X}(d)$. By the multiplication rules in $\mathcal{R}_{n}$, it suffices to check that (\ref{eq:accaoalgmatrizes}) defines an action for basis elements of $\mathcal{R}_{n}$ in the same block $\mathcal{M}_{\binom{n}{r}}(\k S_{r})$. The case $r=0$ is trivial. If $r\geq 1$, let $\sigma, \tau\in S_{r}$ and $I, J, K, L\subseteq \mathbf{n}$ be sets of size $r$. Then
\begin{eqnarray*}
(e_{\alpha}^{\otimes} \centerdot (\sigma E_{I, J}))\centerdot (\tau E_{K, L}) & = & \delta_{X, I} \delta_{K, J} e^{\otimes}_{(\alpha\iota_{I}\sigma \iota_{J}^{-})(\iota_{K}\tau \iota_{L}^{-})} =  \delta_{K, J} (\delta_{X, I} e^{\otimes}_{\alpha \iota_{I}\sigma \tau \iota_{L}^{-}}) \\
& =& \delta_{K, J} (e_{\alpha}^{\otimes} \centerdot (\sigma\tau E_{I, L})) = e_{\alpha}^{\otimes} \centerdot ((\sigma E_{I, J})(\tau E_{K, L})).
\end{eqnarray*}
We also have that $e_{\alpha}^{\otimes} \centerdot 1_{\mathcal{R}_{n}} = e_{\alpha}^{\otimes}\centerdot \epsilon_{\mathbf{r}} E_{X,X} = e^{\otimes}_{\alpha\iota_{X}\epsilon_{\mathbf{r}} \iota_{X}^{-}}=e_{\alpha}^{\otimes}$, where
$$1_{\mathcal{R}_{n}}=\sum_{r=0}^{n} \sum_{|Y|=r} \epsilon_{\mathbf{r}}E_{Y,Y}.$$

This proves that Expression (\ref{eq:accaoalgmatrizes}) defines a right action of $\mathcal{R}_{n}$ on $\otimes^{n} U$.

As to the second statement, let $X\subseteq \mathbf{n}$ with $|X|=r$, $\beta\in \Gamma_{X}(d)$ and $\xi\in\mathcal{S}_{\, \scriptscriptstyle{\k}}(d, \mathbf{n})$. It is enough to prove that $(\xi e_{\beta}^{\otimes})\centerdot (\sigma E_{I,J})= \xi(e_{\beta}^{\otimes}\centerdot (\sigma E_{I,J}))$ for all $\sigma \in S_{r}$ and all $I, J\subseteq \mathbf{n}$ of size $r$. Since, for a fixed $\sigma \in S_{r}$, $\alpha \in \Gamma_{I}(d)$ if and only if $\gamma=\alpha\iota_{I}\sigma \iota_{J}^{-}\in \Gamma_{J}(d)$ and, in such case, $c_{\alpha, \beta}=c_{\gamma, \beta\iota_{I}\sigma\iota_{J}^{-}}$, we have
\begin{eqnarray*}
(\xi e_{\beta}^{\otimes})\centerdot (\sigma E_{I,J}) & = & \sum_{\alpha\in \Gamma_{X}(d)} \xi(c_{\alpha, \beta}) (e_{\alpha}^{\otimes}\centerdot (\sigma E_{I,J}))= \sum_{\alpha\in \Gamma_{X}(d)} \xi(c_{\alpha, \beta})(\delta_{X,I} e_{\alpha\iota_{I}\sigma \iota_{J}^{-}}^{\otimes}) \\
& = &\sum_{\gamma\in \Gamma_{X}(d)}  \delta_{X,J}  \xi(c_{\gamma, \beta\iota_{I}\sigma\iota_{J}^{-}})e_{\gamma}^{\otimes} = \xi(\delta_{X,J} e_{\beta\iota_{I}\sigma \iota_{J}^{-}}^{\otimes}) = \xi(e_{\beta}^{\otimes}\centerdot (\sigma E_{I,J})).
\end{eqnarray*}
\end{proof}

\vspace{0.2cm}

The next two results are a Schur--Weyl duality analog for the extended Schur algebra $\mathcal{S}_{\, \scriptscriptstyle{\k}}(d, \mathbf{n})$ and the matrix algebra $\mathcal{R}_{n} = \displaystyle{\bigoplus_{r=0}^{n}  \mathcal{M}_{\binom{n}{r}}(\k S_{r})}$ on $\otimes^{n} U$.

\vspace{0.2cm}

\begin{theorem}\label{teo:schurweylalgmat}
Let $U$ be an $\k$-space of dimension $d+1$. Let $\mathcal{S}_{\, \scriptscriptstyle{\k}}(d, \mathbf{n})$ act on $\otimes^{n} U$ as in Equation (\ref{eq:accaotensorS}) and let $\rho: (\mathcal{R}_{n})^{\circ} \rightarrow {\rm End}_{\, \k}(\otimes^{n} U)$ be the representation defined in Lemma \ref{le:accaoalgmatrizes}. If $d\geq n$, then 
$$\rho: (\mathcal{R}_{n})^{\circ} \rightarrow {\rm End}_{\, \mathcal{S}_{\, \scriptscriptstyle{\k}}(d, \mathbf{n})}(\otimes^{n} U)$$
is an isomorphism of $\k$-algebras.
\end{theorem}
\begin{proof} We start by showing that $\dim_{\, \scriptscriptstyle{\k}} (\mathcal{R}_{n})^{\circ} = \dim_{\, \scriptscriptstyle{\k}} \left ( {\rm End}_{\, \mathcal{S}_{\, \scriptscriptstyle{\k}}(d, \mathbf{n})}(\otimes^{n} U) \right )$. This follows from Lemma \ref{lem:distributividadeSmod} since, as an $\k$-algebra, ${\rm End}_{\, \mathcal{S}_{\, \scriptscriptstyle{\k}}(d, \mathbf{n})}(\otimes^{n} U)$ is isomorphic to 
\begin{eqnarray*}
{\rm Hom}_{\, \mathcal{S}_{\, \scriptscriptstyle{\k}}(d, \mathbf{n})} \left ( \bigoplus_{r=0}^{n} \binom{n}{r} \otimes^{r} V; \bigoplus_{r=0}^{n} \binom{n}{r} \otimes^{r} V \right ) & \cong & \bigoplus_{r=0}^{n}\mathcal{M}_{\binom{n}{r}}({\rm End}_{\, \mathcal{S}_{\, \scriptscriptstyle{\k}}(d, \mathbf{n})}(\otimes^{r} V)) \\
& \cong & \bigoplus_{r=0}^{n}\mathcal{M}_{\binom{n}{r}}( (\k S_{r})^{\circ}).
\end{eqnarray*}
Thus $\dim_{\, \scriptscriptstyle{\k}} \left ( {\rm End}_{\, \mathcal{S}_{\, \scriptscriptstyle{\k}}(d, \mathbf{n})}(\otimes^{n} U) \right )= \displaystyle{\sum_{r=0}^{n} \binom{n}{r}^{2} r! }= \dim_{\, \scriptscriptstyle{\k}} ( \mathcal{R}_{n} ) = \dim_{\, \scriptscriptstyle{\k}} ( (\mathcal{R}_{n})^{\circ})$.

It remains to show that $\rho: (\mathcal{R}_{n})^{\circ} \rightarrow {\rm End}_{\, \mathcal{S}_{\, \scriptscriptstyle{\k}}(d, \mathbf{n})}(\otimes^{n} U)$ is injective. Let $M\in \mathcal{R}_{n}$ be such that $\rho(M)=0$. By Expression (\ref{eq:basealgebramat}), there are well-determined $a_{I,J}^{\sigma}\in \k$ such that
\begin{equation*}
M= \sum_{k=0}^{n} \sum_{\begin{smallmatrix} I, J\subseteq \mathbf{n}, \\ |I|=|J|=k \end{smallmatrix}} \sum_{\sigma \in S_{k}} a_{I,J}^{\sigma} (\sigma E_{I,J}).
\end{equation*}

Fix $0\le r\leq n$ and choose an arbitrary $X\subseteq \mathbf{n}$ of size $r$. Let $\alpha = \epsilon_{X}\in \Gamma_{X}(d)$. Since $d\geq n$, it follows that $\epsilon_{X}$ is a well-defined element of $\Gamma_{X}(d)$. Then $e_{\alpha}^{\otimes}\in \otimes^{n} U$ and $e_{\alpha}^{\otimes}\centerdot M =0$. This implies that
\begin{equation*}
e^{\otimes}_{\alpha}\centerdot M= \sum_{k=0}^{n} \sum_{\begin{smallmatrix} I, J\subseteq \mathbf{n}, \\ |I|=|J|=k \end{smallmatrix}} \sum_{\sigma \in S_{k}} a_{I,J}^{\sigma} e_{\alpha}^{\otimes}\centerdot (\sigma E_{I,J})= \sum_{\begin{smallmatrix} J\subseteq \mathbf{n}, \\ |J|=r \end{smallmatrix}} \sum_{\sigma \in S_{r}} a_{X,J}^{\sigma} e_{\iota_{X}\sigma\iota_{J}^{-}}^{\otimes} = 0.
\end{equation*}

For any $J\subseteq \mathbf{n}$ of size $r$ and any $\sigma, \tau\in S_{r}$, $\iota_{X}\sigma \iota_{J}^{-}=\iota_{X}\tau \iota_{J}^{-}$ if and only if $\sigma = \tau$. Moreover, for any $K\subseteq \mathbf{n}$ of size $r$, $\gamma\in \Gamma_{K}(d)$ and $\beta\in \Gamma_{J}(d)$, we have that $e_{\gamma}^{\otimes}=e_{\beta}^{\otimes}$ if and only if $K=J$ and $\gamma = \beta$. Hence, the left-hand side of the above equation is a linear combination of distinct elements of the $\k$-basis $\{e_{\nu}^{\otimes} : \nu\in \Gamma_{X}(d), X\subseteq \mathbf{n}\}$ of $\otimes^{n} U$. We deduce that $a_{X, J}^{\sigma}=0$, for all $\sigma \in S_{r}$ and $J\subseteq \mathbf{n}$ of size $r$. Since $r$ and $X$ were chosen arbitrarially, we conclude that $a_{X, J}^{\sigma}=0$, for all $0\leq r\leq n$, all $\sigma \in S_{r}$ and all $X, J\subseteq \mathbf{n}$ such that $|X|=|J|=r$. Hence $M=0$.
\end{proof}

\vspace{0.2cm}

\begin{corollary}\label{cor:schurweylalgmat}
Let $U = V\oplus \k e_{\infty}$ be an $\k$-space of dimension $d+1$ such that $V$ is a $d$-dimensional $\k$-subspace of $U$. The centralizer algebra of $\mathcal{R}_{n}$ on $\otimes^{n} U$ is isomorphic to $\mathcal{S}_{\, \scriptscriptstyle{\k}}(d, \mathbf{n})$, that is, as $\k$-algebras, 
$$\mathcal{S}_{\, \scriptscriptstyle{\k}}(d, \mathbf{n}) \cong {\rm End}_{\, \mathcal{R}_{n}}(\otimes^{n} U).$$
\end{corollary}
\begin{proof} Since $\mathcal{R}_{n}$ is a finite-dimensional (split) semisimple $\k$-algebra, the assertion follows from Theorem \ref{teo:schurweylalgmat} and the Double Centralizer Theorem (Theorem \ref{teo:doublecent}).
\end{proof}

\vspace{0.2cm}

We now make use of the isomorphism of $\k$-algebras of Theorem \ref{teo:rnisomorfoalg} to establish the promised Schur--Weyl duality between $\mathcal{S}_{\, \scriptscriptstyle{\k}}(d, \mathbf{n})$ and $\k R_{n}$ on $\otimes^{n} U$.

\vspace{0.2cm}

\begin{theorem} \label{teo:ligacaoSRn}
Let $U$ be an $\k$-vector space of dimension $d+1$ and let $\mathcal{S}_{\, \scriptscriptstyle{\k}}(d, \mathbf{n})$ act on $\otimes^{n} U$ as in Equation (\ref{eq:accaotensorS}). For all $\sigma \in R_{n}$ and all $\alpha\in \Gamma_{X}(d)$ with $X\subseteq \mathbf{n}$, define
\begin{equation}
\label{eq:accaoRneltosbase}
e_{\alpha}^{\otimes} \cdot \sigma = 
 \begin{cases}
 e_{\alpha\sigma}^{\otimes} & \textrm{if $X \subseteq R(\sigma)$} \\ 
0 & \textrm{otherwise.}
\end{cases}
\end{equation}
If $d\geq n$, then the left action of $\mathcal{S}_{\, \scriptscriptstyle{\k}}(d, \mathbf{n})$ and the right action of $\k R_{n}$ defined by Equation (\ref{eq:accaoRneltosbase}) on $\otimes^{n} U$ generate the full centralizers of each other on ${\rm End}_{\, \k}(\otimes^{n} U)$. In addition, there are isomorphisms of $\k$-algebras such that
$$\mathcal{S}_{\, \scriptscriptstyle{\k}}(d, \mathbf{n}) \cong {\rm End}_{\, \k R_{n}} (\otimes^{n} U) \ {\rm \ and \ } \ (\k R_{n})^{\circ} \cong {\rm End}_{\, \mathcal{S}_{\, \scriptscriptstyle{\k}}(d, \mathbf{n})} (\otimes^{n} U).$$
\end{theorem}
\begin{proof}
Let $\phi: \k R_{n} \rightarrow \mathcal{R}_{n}$ be the $\k$-algebra isomorphism of Theorem \ref{teo:rnisomorfoalg}. By Lemma \ref{le:accaoalgmatrizes} and Theorem \ref{teo:rnisomorfoalg}, the rule $z \cdot \sigma = z \centerdot \phi(\sigma)$, for all $z\in \otimes^{n} U$ and $\sigma \in R_{n}$, turns $\otimes^{n} U$ into a right $\k R_{n}$-module. With this right $\k R_{n}$-action, the result follows from Theorem \ref{teo:schurweylalgmat} and Corollary \ref{cor:schurweylalgmat} and the fact that $\phi$ is an isomorphism of $\k$-algebras. The previous right $\k R_{n}$-action reduces to Expression (\ref{eq:accaoRneltosbase}) for basis elements of $e_{\alpha}^{\otimes}$ of $\otimes^{n} U$, where $\alpha \in\Gamma_{X}(d)$ for some $X\subseteq \mathbf{n}$. This can be verified by $$e_{\alpha}^{\otimes} \cdot \sigma = e_{\alpha}^{\otimes} \centerdot \phi(\sigma)  = \sum_{J\subseteq D(\sigma)} e_{\alpha}^{\otimes} \centerdot (\mathfrak{p}(\sigma\epsilon_{J}) E_{\sigma(J), J}) =  \sum_{J\subseteq D(\sigma)} \delta_{X,J} e^{\otimes}_{\alpha \epsilon_{\sigma(J)} \sigma \epsilon_{J}}.$$ 
\end{proof}

\vspace{0.2cm}

Let's illustrate how an element of $R_{n}$ acts on $\otimes^{n} U$ with an example. If $d=n=5$ and $\alpha=(5,2,2)\in \Gamma_{X}(5)$ with $X=\{1,4,5\}$, the corresponding basis element of $\otimes^{5} U$ is given by
$$e_{\alpha}^{\otimes}=e_{5}\otimes e_{\infty} \otimes e_{\infty} \otimes e_{2} \otimes e_{2}.$$

If $\sigma \in R_{5}$ is the element $\sigma =
\left ( \begin{array}{ccccc}
1 & 2 & 3 & 4 & 5 \\
5 & - & 1 & 2 & 4
\end{array} \right )$, then $X\subseteq R(\sigma)=\{1,2,4,5\}$ and $\alpha \sigma = (2,5,2)\in \Gamma_{Y}(5)$, where $Y=\{1,3,5\}$. Hence,
$$e_{\alpha}^{\otimes} \cdot \sigma = e_{\alpha\sigma}^{\otimes} = e_{2} \otimes e_{\infty} \otimes e_{5} \otimes e_{\infty} \otimes e_{2} \in \otimes^{5} U.$$

On the other hand, if $\sigma \in S_{n}\subseteq R_{n}$, Expression (\ref{eq:accaoRneltosbase}) becomes $e_{\alpha}^{\otimes} \cdot \sigma = e_{\alpha\sigma}^{\otimes}$, for all $\alpha\in \Gamma_{X}(d)$ and $X\subseteq \mathbf{n}$. This is the usual right $S_{n}$-action by place permutations on $\otimes^{n} U$. Indeed, if $u_{1}, \ldots , u_{n}$ are vectors in $U$, there are well-determined $u_{\alpha, X}\in \k$ such that $u_{1}\otimes \ldots \otimes u_{n} = \displaystyle{\sum_{X\subseteq \mathbf{n}} \sum _{\alpha \in \Gamma_{X}(d)} u_{\alpha, X} e_{\alpha}^{\otimes}}$. Since $\sigma \in S_{n}$, we have that 
$$u_{\sigma(1)}\otimes \ldots \otimes u_{\sigma(n)} = \sum_{Y\subseteq \mathbf{n}} \sum_{\beta \in \Gamma_{Y}(d)} u_{\beta\sigma^{-1}, \sigma(Y)} e_{\beta}^{\otimes} = \sum_{X\subseteq \mathbf{n}} \sum_{\alpha \in \Gamma_{X}(d)} u_{\alpha, X} e_{\alpha\sigma}^{\otimes},$$
an equality obtained by reindexation ($\alpha=\beta\sigma^{-1} \Leftrightarrow \beta=\alpha\sigma$ and $\sigma(Y)=X$, for all $X\subseteq \mathbf{n}$). This implies that 
$$(u_{1}\otimes \ldots \otimes u_{n} ) \cdot \sigma = u_{\sigma(1)}\otimes \ldots \otimes u_{\sigma(n)}.$$
Thus, the restriction of the right $R_{n}$-action on $\otimes^{n} U$ given by Expression (\ref{eq:accaoRneltosbase}) to $S_{n}$ is the usual right action of $S_{n}$ on $\otimes^{n} U$ by place permutations.

Although it is a lengthy computation, it is possible to show that the $R_{n}$-action on $\otimes^{n} U$ by ``place permutations" defined by L. Solomon in \cite[Equation $(5.5)$]{Solomon} turns into Expression (\ref{eq:accaoRneltosbase}) for basis elements of $\otimes^{n} U$. The reader should be aware that Solomon makes use of a different convention than ours when it comes to compose in $R_{n}$. This means that Theorem \ref{teo:ligacaoSRn} is a reformulation in the setting of Schur algebras of Solomon's important Schur--Weyl duality between $G_{d}$ and $R_{n}$ on $\otimes^{n} U$ \cite[Theorem $5.10$]{Solomon}.

It is worth noting that our approach opens up new possibilities for better understanding the extent of the interactions between the representation theories of rook monoids, general linear groups, symmetric groups and (extended) Schur algebras.

On one hand, our starting point was to view $G_{d}$ as a subgroup of $G_{d+1}$. This implies that a simple $\k G_{d+1}$-module is also an $\k G_{d}$-module. Hence, it makes sense considering its decomposition into simple $\k G_{d}$-modules. In the language of Schur algebras, this amounts to decomposing a simple $\mathcal{S}_{\, \scriptscriptstyle{\k}}(d+1, n)$-module into simple $\mathcal{S}_{\, \scriptscriptstyle{\k}}(d, \mathbf{n})$-modules and determining the corresponding multiplicities. This procedure is known as a {\it branching rule} for $\mathcal{S}_{\, \scriptscriptstyle{\k}}(d, \mathbf{n}) \subseteq \mathcal{S}_{\, \scriptscriptstyle{\k}}(d+1, n)$ (see \cite{Halverson96} and \cite[Chapter $8$]{GoodmanWallach}). Since both $\mathcal{S}_{\, \scriptscriptstyle{\k}}(d, \mathbf{n})$ and $\mathcal{S}_{\, \scriptscriptstyle{\k}}(d+1, n)$ are finite-dimensional semisimple $\k$-algebras, the branching rule for $\mathcal{S}_{\, \scriptscriptstyle{\k}}(d, \mathbf{n}) \subseteq \mathcal{S}_{\, \scriptscriptstyle{\k}}(d+1, n)$ is the same as that for ${\rm End}_{\, \mathcal{S}_{\, \scriptscriptstyle{\k}}(d+1, n)} \left ( \otimes^{n} U \right)  \subseteq {\rm End}_{\,\mathcal{S}_{\, \scriptscriptstyle{\k}}(d, \mathbf{n})} \left ( \otimes^{n} U \right )$ \cite[Theorem $1.7.$]{Halverson96}. By Theorem \ref{teo:ligacaoSRn}, this means that it is possible to derive concise proofs of combinatorial formulas for multiplicities associated with the restriction to $S_{n} \subseteq R_{n}$ of irreducible characters of $R_{n}$. In the near future, we hope to publish these proofs and recover in an economical way some of the results in \cite[Section $3$]{Solomon}.

Another upshot of our approach is that we may use the tools associated with Schur algebras to give a construction of the irreducible modules of the rook monoid realized on tensors which is analogous to that of the dual Specht modules for the symmetric group. Indeed, we may apply Green's techniques \cite[Chapter $6$]{Green1} and define an idempotent $\zeta \in \mathcal{S}_{\, \scriptscriptstyle{\k}}(d, \mathbf{n})$ which satisfies the algebra isomorphism $\zeta \mathcal{S}_{\, \scriptscriptstyle{\k}}(d, \mathbf{n}) \zeta \cong \k R_{n} $. The idempotent $\zeta$ induces a functor between the module categories for $\mathcal{S}_{\, \scriptscriptstyle{\k}}(d, \mathbf{n})$ and $\k R_{n}$ and we can make use of this functor to build a complete set of simple modules for $\k R_{n}$ from the Carter-Lusztig $\mathcal{S}_{\, \scriptscriptstyle{\k}}(d, \mathbf{n})$-modules \cite{CarterLusztig} realized on $\otimes^{n} U$, obtaining an analog for $R_{n}$ of the dual Specht modules for $S_{n}$ \cite{Green1}. We also expect to exhibit this construction in the near future.

Finally, in this article, we have laid the foundations for studying the modular representations of the rook monoid on tensors. Although Theorem \ref{teo:rnisomorfoalg} was stated in characteristic zero and our Schur-Weyl duality between $R_{n}$ and $\mathcal{S}_{\, \scriptscriptstyle{\k}}(d, \mathbf{n})$ on $\otimes^{n} U$ relies heavily on this result, it is possible to show that a Schur--Weyl duality between (a subalgebra of) $\k R_{n}$ and $\mathcal{S}_{\, \scriptstyle{\k}}(d, \mathbf{n})$ on a tensor space can be established for infinite fields.


\vspace{0.3cm}

\bibliography{bibliography}

\end{document}